\documentclass[11pt]{amsart}
\usepackage[latin1]{inputenc}

\usepackage{bbm,bm}
\usepackage{amssymb}
\usepackage{stmaryrd}

\usepackage{subfigure}
\usepackage[width=.85\textwidth,font=small,labelfont=sc]{caption}

\PII{}
\copyrightinfo{}{}

\usepackage{color}
\usepackage{pstricks}

\usepackage{xspace}
\usepackage{enumitem}

\usepackage[bookmarks=false,colorlinks,linkcolor=mylinkcolor,citecolor=mycitecolor]{hyperref}
	\definecolor{mycitecolor}{rgb}{.1,.6,.1}
	\definecolor{mylinkcolor}{rgb}{.1,.1,.6}

\newtheorem{app}{Application}

\newtheorem{corollary}[app]{Corollary}
\newtheorem{lemma}[app]{Lemma}
\newtheorem{proposition}[app]{Proposition}
\newtheorem{theorem}[app]{Theorem}

\newtheorem{remark}[app]{Remark} 

\newcommand{\hilb}[2]{%
   \mathop{H_{#1}}\bigl(#2\bigr)
}

\newcommand{\lexl}{\mathrel{<_{\mathrm{lex}}}}

\newcommand{\lexgeq}{\mathrel{\geq_{\mathrm{lex}}}}
\newcommand{\PI}{\mathsf{B}}
\newcommand{\PPI}{\mathsf{PB}}
\newcommand{\Comp}{\mathsf{C}}
\newcommand{\sS}{\mathfrak{C}}

\newcommand{\lrct}{Littlewood--Richardson composition tableau\xspace}
\newcommand{\lrctx}{Littlewood--Richardson composition tableaux\xspace}
\newcommand{\id}{\mathbbm1}
\newcommand{\shape}{\bm\lambda}
\newcommand{\length}{\bm{l}}
\newcommand{\std}{\mathrm{st}}
\newcommand{\destd}{\mathbf{d}}
\newcommand{\Destd}{\mathsf{D}}
\newcommand{\rdes}{\mathbf{r}}
\newcommand{\HilbQuo}{P}
\newcommand{\bN}{\mathbb{N}}

\newcommand{\bQ}{\mathbb{Q}}

\newcommand{\bZ}{\mathbb{Z}}

\newcommand{\cB}{\mathcal{B}}

\newcommand{\cE}{\mathcal{E}}

\newcommand{\cS}{\mathcal{S}}

\definecolor{gr90}{gray}{0.90}

\definecolor{gr75}{gray}{0.75}

\definecolor{gyblue}{cmyk}{0,0.5,0,0}

\newcommand{\Sym}{\mathit{Sym}}
\newcommand{\Qsym}{\mathit{QSym}}

\newcommand{\frakS}{\mathfrak{S}}
\newcommand{\xx}{\mathbf{x}}

\newcommand{\suchthat}{:}

\newcommand{\qschur}{\mathcal{S}}

\DeclareMathOperator{\sM}{\mathfrak{M}}

\newlength{\cellsize} 
\newsavebox{\cell}

\sbox{\cell}{%
\begin{picture}(18,18)\linethickness{0.6pt} %
  \put(0,0){\line(1,0){18}} \put(0,0){\line(0,1){18}}
  \put(18,0){\line(0,1){18}} \put(0,18){\line(1,0){18}}
\end{picture}}

\newcommand\cellify[1]{%
  \def\thearg{#1}\def\nothing{}%
  \ifx\thearg\nothing \vrule width0pt height\cellsize depth0pt\else
  \hbox to 0pt{\usebox{\cell} \hss}\fi%
  \vbox to \cellsize{ \vss \hbox to
  \cellsize{\hss$\scriptstyle#1$\hss} \vss}
}
\newcommand\tableau[1]{\vtop{\let\\\cr
\baselineskip -16000pt \lineskiplimit 16000pt \lineskip 0pt
\ialign{&\cellify{##}\cr#1\crcr}}}

\newcommand\bas[1]{\omit \vbox to \cellsize{ \vss \hbox to \cellsize{\hss$#1$\hss} \vss}}

\author{Aaron Lauve}
\address[Lauve]{
	Department of Mathematics\\
        Texas A\&M University\\
        College Station, TX \ 77843
        }
\email{lauve@math.tamu.edu}
\urladdr{http://www.math.tamu.edu/\~{}lauve}

\author{Sarah K Mason}
\address[Mason]{
	Department of Mathematics\\
        Wake Forest University\\
        Winston-Salem, NC \ 27109
	}
\email{sarahkmason@gmail.com}
\urladdr{http://sarahmason.wikidot.com/home}
\thanks{Mason was partially supported by NSF postdoctoral fellowship DMS 0603351.}

\title{QSym over Sym has a stable basis}


\begin{document}

\begin{abstract}
We prove that the subset of quasisymmetric polynomials conjectured by Bergeron and Reutenauer to be a basis for the coinvariant space of quasisymmetric polynomials is indeed a basis.  This provides the first constructive proof of the Garsia--Wallach result stating that quasisymmetric polynomials form a free module over symmetric polynomials and that the dimension of this module is $n!$.
\end{abstract}

\maketitle

\section{Introduction}\label{sec: intro}

Quasisymmetric polynomials have held a special place in algebraic combinatorics since their introduction in \cite{Ges:1984}.  They are the natural setting for many enumeration problems \cite{Sta:1999} as well as the development of Dehn--Somerville relations \cite{AguBerSot:2006}.  In addition, they are related in a natural way to Solomon's descent algebra of the symmetric group \cite{MalReu:1995}. In this paper, we follow \cite[Chapter 11]{Ber:2009} and view them through the lens of invariant theory. Specifically, we consider the relationship between the two subrings $\Sym_n\subseteq \Qsym_n\subseteq \bQ[\xx]$ of symmetric and quasisymmetric polynomials in variables $\xx=\xx_n:=\{x_1,x_2,\dotsc,x_n\}$. Let $\cE_n$ denote the ideal in $\Qsym_n$ generated by the elementary symmetric polynomials. In 2002, F. Bergeron and C. Reutenauer made a sequence of three successively finer conjectures concerning the quotient ring $\Qsym_n/\cE_n$. A. Garsia and N. Wallach were able to prove the first two in \cite{GarWal:2003}, but the third one remained open; we close it here (Corollary \ref{thm: m-basis}) with the help of a new basis for $\Qsym_n$ introduced in \cite{HLMvW:1}.

\subsection*{Acknowledgements} We thank Fran\c{c}ois Bergeron for sharing the preceding story with us and encouraging us to write this paper. Our approach follows an idea that he proposed during CanaDAM 2009. We also gratefully acknowledge several beneficial discussions we had with Adriano Garsia, Christophe Reutenauer, and Frank Sottile.

\subsection{Motivating context}
Recall that $\Sym_n$ is the ring $\bQ[\xx]^{\frakS_n}$ of invariant polynomials under the permutation action of $\frakS_n$ on $\xx$ and $\bQ[\xx]$. One of the crowning results in the invariant theory of $\frakS_n$ is that the following true statements are equivalent:
\begin{enumerate}[label={\small(S\arabic*)}, ref={\small S\arabic*}]
\item \label{itm: poly}
$\bQ[\xx]^{\frakS_n}$ is a polynomial ring, generated, say, by the elementary symmetric polynomials
$\cE_n=\{e_1(\xx), \dotsc, e_n(\xx)\}$;
\item \label{itm: free-module}
the ring $\bQ[\xx]$ is a free $\bQ[\xx]^{\frakS_n}$-module;
\item \label{itm: coinvariant}
the \emph{coinvariant space} $\bQ[\xx]_{\frakS_n} = \bQ[\xx] /\bigl(\cE_n\bigr)$
has dimension $n!$ and is isomorphic to the regular representation of $\frakS_n$.
\end{enumerate}
See \cite[\S\S17, 18]{Kan:2001} for details. Analogous statements hold on replacing $\frakS_n$ by any pseudo-reflection group. Since all spaces in question are graded,
we may add a fourth item to the list: the \emph{Hilbert series}
$\hilb{q}{\bQ[\xx]_{\frakS_n}} = \sum_{k\geq0} d_k\,q^k$,
where $d_k$ records the dimension of the $k$th graded component of $\bQ[\xx]_{\frakS_n}$, satisfies
%
\begin{enumerate}[resume*]
\item \label{itm: hilb-quotient}
$\displaystyle \hilb{q}{\bQ[\xx]_{\frakS_n}} = \hilb{q}{\bQ[\xx]} \Big/ \hilb{q}{\bQ[\xx]^{\frakS_n}} .$
\end{enumerate}

Before we formulate the conjectures of Bergeron and Reutenauer, we recall another page in the story of $\Sym_n$ and the quotient space $\bQ[\xx]/\left(\cE_n\right)$. The ring homomorphism $\zeta$ from $\bQ[\xx_{n+1}]$ to $\bQ[\xx_n]$ induced by the mapping $x_{n+1}\mapsto 0$ respects the rings of invariants (that is, $\zeta\colon \Sym_{n+1} \twoheadrightarrow \Sym_n$ is a ring homomorphism). Moreover, $\zeta$ respects the fundamental bases of monomial ($m_\lambda$) and Schur ($s_\lambda$) symmetric polynomials of $\Sym_n$, indexed by partitions $\lambda$ with at most $n$ parts. For example,
\[
  \zeta(m_\lambda(\xx_{n+1})) = \begin{cases}
    m_\lambda(\xx_n), & \hbox{if }\lambda\hbox{ has at most $n$ parts}, \\
    0, & \hbox{otherwise.}
   \end{cases}
\]
The stability of these bases plays a crucial role in representation theory \cite{Mac:1995}. Likewise, the associated stability of bases for the coinvariant spaces (e.g., of Schubert polynomials \cite{FomKir:1996,LasSch:1982,Man:2001}) plays a role in the cohomology theory of flag varieties.

\subsection{Bergeron--Reutenauer context}
Given that $\Qsym_n$ is a polynomial ring \cite{MalReu:1995} containing $\Sym_n$, one might ask, by analogy with $\bQ[\xx]$, how $\Qsym_n$ looks as a module over $\Sym_n$. This was the question investigated by Bergeron and Reutenauer \cite{BerReu}. They began by computing the quotient $\HilbQuo_n(q) := \hilb{q}{\Qsym_n} \big/ \hilb{q}{\Sym_n}$ by analogy with \eqref{itm: hilb-quotient}. Surprisingly, the result was a polynomial in $q$ with nonnegative integer coefficients (so it could, conceivably, enumerate the graded space $\Qsym_n/\cE_n$). More astonishingly, sending $q$ to $1$ gave $\HilbQuo_n(1) = n!$. This led to the following two conjectures, subsequently proven in \cite{GarWal:2003}:
\begin{enumerate}[label={\small(Q\arabic*)}, ref={\small Q\arabic*}]
\item \label{conj:free}
  The ring $\Qsym_n$ is a free module over $\Sym_n$;
\item \label{conj:dimension}
  The dimension of the ``coinvariant space'' $\Qsym_n/\cE_n$ is $n!$.
\end{enumerate}
%

In their efforts to solve the conjectures above, Bergeron and Reutenauer introduced the notion of ``pure and inverting'' compositions $\PI_n$ with at most $n$ parts. These compositions have the favorable property of being $n$-stable in that $\PI_n \subseteq \PI_{n+1}$ and that $\PI_{n+1} \setminus \PI_n$ are the pure and inverting compositions with exactly $n{+}1$ parts. They were able to show that the pure and inverting ``quasi-monomials'' $M_\beta$ (see Section \ref{sec: qsym}) span $\Qsym_n/\cE_n$ and that they are $n!$ in number. However, the linear independence of these polynomials over $\Sym_n$ remained open. Their final conjecture, which we prove in Corollary \ref{thm: m-basis}, is as follows:
\begin{enumerate}[resume*]
\item \label{conj:basis}
  The set of quasi-monomials $\left\{ M_\beta \suchthat \beta \in \PI_n \right\}$ is a basis for $\Qsym_n/\cE_n$.
\end{enumerate}

The balance of this paper is organized as follows. In Section \ref{sec: qsym}, we recount the details surrounding a new basis $\{\cS_\alpha\}$ for $\Qsym_n$ called the quasisymmetric Schur polynomials. These behave particularly well with respect to the $\Sym_n$ action in the Schur basis. In Section \ref{sec: coinvariant-space}, we give further details surrounding the ``coinvariant space'' $\Qsym_n / \cE_n$. These include a bijection between compositions $\alpha$ and pairs $(\lambda,\beta)$, with $\lambda$ a partition and $\beta$ a pure and inverting composition, that informs our main results. Section \ref{sec: main-results} contains these results---a proof of \eqref{conj:basis}, but with the quasi-monomials $M_\beta$ replaced by the quasisymmetric Schur polynomials $\cS_\beta$. We conclude in Section \ref{sec: corollaries} with some corollaries to the proof. These include \eqref{conj:basis} as originally stated, as well as a version of \eqref{conj:free} and \eqref{conj:basis} over the integers.

\section{Quasisymmetric polynomials}\label{sec: qsym}

A polynomial in $n$ variables $\xx=\{x_1,x_2, \dotsc, x_n \}$ is said to be {\it quasisymmetric} if and only if for each composition $(\alpha_1, \alpha_2, \dotsc , \alpha_k)$, the monomial $x_1^{a_1} x_2^{a_2} \dotsb x_k^{a_k}$ has the same coefficient as $x_{i_1}^{\alpha_1} x_{i_2}^{\alpha_2} \dotsb x_{i_k}^{\alpha_k}$ for all sequences $1 \le i_1 < i_2 < \dotsb < i_k \le n$.  For example, $x_1^2 x_2 + x_1^2 x_3 + x_2^2 x_3$ is a quasisymmetric polynomial in the variables $\{x_1,x_2,x_3 \}$. The ring of quasisymmetric polynomials in $n$ variables is denoted $QSym_n$. (Note that every symmetric polynomial is quasisymmetric.)

It is easy to see that $\Qsym_n$ has a vector space basis given by the quasi-monomials
\[
   M_\alpha(\xx) = \sum_{i_1<\dotsb < i_k} x_{i_1}^{\alpha_1}\dotsb x_{i_k}^{\alpha_k} ,
\]
for $\alpha=(\alpha_1,\dotsc,\alpha_k)$ running over all compositions with at most $n$ parts. It is less evident that $\Qsym_n$ is a ring, but see \cite{Haz:2003} for a formula for the product of two quasi-monomials.
We write $\length(\alpha)=k$ for the \emph{length} (number of parts) of $\alpha$ in what follows.
%
We return to the quasi-monomial basis in Section \ref{sec: corollaries}, but for the majority of the paper, we focus on the basis of ``quasisymmetric Schur polynomials'' as its known multiplicative properties assist in our proofs.

\subsection{The basis of quasisymmetric Schur polynomials}\label{sec: qschur}

\begin{figure}[hbt]
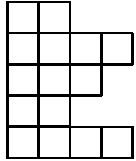

\[
  \tableau{{} & {} \\ {} & {} & {} & {} \\ {} & {} & {} \\
  {} & {} \\ {} & {} & {} & {} }
\]
\caption{The diagram associated to the composition $(2,4,3,2,4)$}
\label{compdiag}
\end{figure}

A {\it quasisymmetric Schur polynomial} $\cS_\alpha$ is defined combinatorially through fillings of composition diagrams. Given a composition $\alpha=(\alpha_1, \alpha_2, \dotsc, \alpha_k)$, its associated diagram is constructed by placing $\alpha_i$ boxes, or {\it cells}, in the $i^{\textit{th}}$ row from the top.  (See Figure \ref{compdiag}.)  The cells are labeled using matrix notation; that is, the cell in the $j^{\textit{th}}$ column of the $i^{\textit{th}}$ row of the diagram is denoted $(i,j)$.  We abuse notation by writing $\alpha$ to refer to the diagram for $\alpha$.

Given a composition diagram $\alpha = (\alpha _1, \alpha _2, \ldots , \alpha _\ell)$ with largest part $m$, a \emph{composition tableau} $T$ of \emph{shape} $\alpha$ is a filling of the cells $(i,j)$ of $\alpha$ with positive integers $T(i,j)$ such that
\begin{enumerate}[label={\small(CT\arabic*)}, ref={\small CT\arabic*}]
\item \label{itm: CT-row} entries in the rows of $T$ weakly decrease when read from left to right,
\item \label{itm: CT-col} entries in the leftmost column of $T$ strictly increase when read from top to bottom,
\item \label{itm: CT-triple} entries satisfy the \emph{triple rule:}

\noindent Let $(i,k)$ and $(j,k)$ be two cells in the same column so that $i < j$.  If $\alpha_i \ge \alpha_j$ then either $T(j,k) < T(i,k)$ or $T(i,k-1) < T(j,k)$.  If $\alpha_i < \alpha_j$ then either $T(j,k) < T(i,k)$ or $T(i,k)<T(j,k+1)$.
\end{enumerate}

Assign a {\it weight, $x^T$} to each composition tableau $T$ by letting $a_i$ be the number of times $i$ appears in $T$ and setting $x^T=\prod x_i^{a_i}$.  The quasisymmetric Schur polynomial $\qschur_{\alpha}$ corresponding to the composition $\alpha$ is defined by
\begin{equation*}
  \qschur_{\alpha}(\xx_n) = \sum_T x^T ,
\end{equation*}
the sum being taken over all composition tableaux $T$ of shape $\alpha$ with entries chosen from $[n]$. (See Figure \ref{fig: qschur}.) Each polynomial $\qschur_{\alpha}$ is quasisymmetric and the collection $\{\qschur_{\alpha} \suchthat \length(\alpha) \le n \}$ forms a basis for $QSym_n$~\cite{HLMvW:1}.
\begin{figure}[hbt]
\[
   \tableau{1 & 1 & 1 \\ 2 \\ 3 & 3} \quad \tableau{1 & 1 & 1 \\ 2 \\ 4 & 3} \quad
   \tableau{1 & 1 & 1 \\ 2 \\ 4 & 4} \quad \tableau{1 & 1 & 1 \\ 3 \\ 4 & 4} \; \; \;
   \tableau{2 & 1 & 1 \\ 3 \\ 4 & 4} \quad \tableau{2 & 2 & 1 \\ 3 \\ 4 & 4} \quad
   \tableau{2 & 2 & 2 \\ 3 \\ 4 & 4}
\]
\caption{The composition tableaux encoded in the polynomial  $\qschur_{(3,1,2)}(\xx_4)=x_1^3 x_2 x_3^2 + x_1^3x_2x_3x_4 + x_1^3x_2x_4^2+x_1^3x_3x_4^2 + x_1^2 x_2 x_3 x_4^2 + x_1x_2^2x_3x_4^2 + x_2^3x_3x_4^2$.}
\label{fig: qschur}
\end{figure}

\subsection{Sym action in the Quasisymmetric Schur polynomial basis}

We need several definitions in order to describe the multiplication rule for quasisymmetric Schur polynomials found in \cite{HLMvW:2}.  First, given two compositions $\alpha=(\alpha_1,\dotsc,\alpha_r)$ and $\beta=(\beta_1,\dotsc,\beta_s)$, we say $\alpha$ \emph{contains} $\beta$ ($\alpha \supseteq \beta$) if $r\geq s$ and there is a subsequence $i_1>\dotsb>i_s$ satisfying $\alpha_{i_1} \geq \beta_1, \dotsc, \alpha_{i_s} \geq \beta_s$. The reverse of a partition $\lambda$ is the composition $\lambda^*$ obtained by reversing the order of its parts.  Symbolically, if $\lambda=(\lambda_1, \lambda_2, \hdots , \lambda_k)$ then $\lambda^*=(\lambda_k , \hdots , \lambda_2 , \lambda_1)$.  Let $\beta$ be a composition, let $\lambda$ be a partition, and let $\alpha$ be a composition obtained by adding $|\lambda|$ cells to $\beta$, possibly between adjacent rows of $\beta$.  A filling of the cells of $\alpha$ is called a {\it Littlewood--Richardson composition tableau} of shape $\alpha\supseteq\beta$ if it satisfies the following rules:

\begin{enumerate}[label={\small(LR\arabic*)}, ref={\small LR\arabic*}]
\item \label{itm: LR-row} The $i^{\textit{th}}$ row from the bottom of $\beta$ is filled with the entries $k+i$.
\item \label{itm: LR-content} The content of the appended cells is $\lambda^*$.
\item \label{itm: LR-CTs} The filling satisfies conditions \eqref{itm: CT-row} and \eqref{itm: CT-triple} from Section \ref{sec: qschur}.
\item \label{itm: LR-lattice} The entries in the appended cells, when read from top to bottom, column by column, from right to left, form a \emph{reverse lattice word.} That is, one for which each prefix contains at least as many $i$'s as $(i-1)$'s for each $1<i\le k$.
\end{enumerate}

The following theorem provides a method for multiplying an arbitrary quasisymmetric Schur polynomial by an arbitrary Schur polynomial.

\begin{proposition}[\cite{HLMvW:2}]\label{thm:LR-qschur}
In the expansion
\begin{equation} \label{eq:LR-qschur}
 s_\lambda(\xx) \cdot \qschur_\alpha(\xx) = \sum_\gamma C^ \gamma_{\lambda\,\alpha} \, \qschur_\gamma(\xx),
\end{equation}
the coefficient $C^ \gamma_{\lambda\,\alpha}$ is the number of \lrctx of
shape $\gamma\supseteq\alpha$ with appended content $\lambda ^{*}$.
\end{proposition}

\section{The coinvariant space for quasisymmetric polynomials}\label{sec: coinvariant-space}

Let $B\subseteq A$ be two $\bQ$-algebras with $A$ a free left module over $B$. This implies the existence of a subset $C\subseteq A$ with $A \simeq B\otimes C$ as vector spaces over $\bQ$. In the classical setting of invariant theory (where $B$ is the subring of invariants for some group action on $A$), this set $C$ is identified as coset representatives for the quotient $A/(B_+)$, where $(B_+)$ is the ideal in $A$ generated by the positive part of the graded algebra $B=\bigoplus_{k\geq0} B_{k}$.

Now suppose that $A$ and $B$ are graded rings. If $A$ is free over $B$, then the Hilbert series of $C$ is given as the quotient $\hilb{q}{A} \big/ \hilb{q}{B}$. Let us try this with the choice $A=\Qsym_n$ and $B=\Sym_n$. It is well-known that the Hilbert series for $\Qsym_n$ and $\Sym_n$ are given by
\begin{align}
\label{eq: qsym}
   \hilb{q}{\Qsym_n} &= 1 + \frac{q}{1-q} + \dotsb + \frac{q^n}{(1-q)^n}  \\
\intertext{and}
\label{eq: sym}
   \hilb{q}{\Sym_n} &= \prod_{i=1}^n \frac1{1-q^i} \,.
\end{align}
Let $\HilbQuo_n(q)=\sum_{k\geq0} p_k\, q^k$ denote the quotient of \eqref{eq: qsym} by \eqref{eq: sym}. It is easy to see that
\[
   \HilbQuo_n(q) = \prod_{i=1}^{n-1}\bigl(1+q+\dotsb+q^i\bigr) \sum_{i=0}^n q^i (1-q)^{n-i} \,,
\]
and hence $\HilbQuo_n(1)=n!$. It is only slightly more difficult (see (0.13) in \cite{GarWal:2003}) to show that $\HilbQuo_n(q)$ satisfies the recurrence relation
\begin{equation}\label{eq: recursion}
   \HilbQuo_n(q) = \HilbQuo_{n-1}(q) + q^n\bigl([n]_q! - \HilbQuo_{n-1}(q)\bigr) ,
\end{equation}
where $[n]_q!$ is the standard $q$-version of $n!$. Bergeron and Reutenauer use this recurrence to show that $p_k$ is a nonnegative integer for all $k\ge0$ and to produce a set of compositions $\PI_n$ satisfying $p_k = \#\{ \beta\in\PI_n \suchthat |\beta|=k\}$ for all $n$. In particular, $|\PI_n|=n!$.

Let $\cE_n$ be the ideal in $\Qsym_n$ generated by all symmetric polynomials with zero constant term and call $R_n:=QSym_n / \cE_n$ the {\it coinvariant space for quasisymmetric polynomials}.
From the above discussion, $R_n$ has dimension at most $n!$. If the set of quasi-monomials $\{M_\beta \in \Qsym_n \suchthat \beta \in \PI_n\}$ are linearly independent over $\Sym_n$, then it has dimension exactly $n!$ and $\Qsym_n$ becomes a free $\Sym_n$ module of the same dimension.

\subsection{Destandardization of permutations}

To produce a set $\PI_n$ of compositions indexing a proposed basis of $R_n$, first recognize the $[n]_q!$ in \eqref{eq: recursion} as the Hilbert series for the classical coinvariant space $\bQ[\xx]\big/(\cE_n)$ from \eqref{itm: coinvariant}. The standard set of compositions indexing this space are the \emph{Artin monomials} $\{x_1^{\alpha_1}\dotsb x_n^{\alpha_n} \suchthat 0\leq \alpha_i\leq n-i \}$, but these do not fit into the desired recurrence \eqref{eq: recursion} with $n$-stability. In \cite{Gar:1980}, Garsia developed an alternative set of monomials indexed by permutations. His ``descent monomials'' (actually, the ``reversed'' descent monomials, see \cite[\S 6]{GarWal:2003}) were chosen as the starting point for the recursive construction of the sets $\PI_n$. Here we give a description in terms of ``destandardized permutations.''

In what follows, we view partitions and compositions as words in the alphabet $\bN=\{0,1,2,\dotsc\}$. For example, we write $2543$ for the composition $(2,5,4,3)$. The \emph{standardization} $\std(w)$ of a word $w$ of length $k$ is a permutation in $\frakS_k$ obtained by first replacing (from left to right) the $\ell_1$ $1$s in $w$ with the numbers $1,\dotsc,\ell_1$, then replacing (from left to right) the $\ell_2$ $2$s in $w$ with the numbers $\ell_1{+}1,\dotsc,\ell_1{+}\ell_2$, and so on. For example, $\std(121)=132$ and $\std(2543) = 1432$. The \emph{destandardization} $\destd(\sigma)$ of a permutation $\sigma\in\frakS_k$ is the lexicographically least word $w \in (\bN_+)^k$ satisfying $\std(w)=\sigma$. For example, $\destd(132) = 121$ and $\destd(1432) = 1321$. Let $\Destd_{(n)}$ denote the compositions $\{\destd(\sigma) \suchthat \sigma\in\frakS_n \}$. Finally, given $\destd(\sigma)=(\alpha_1,\dotsc,\alpha_k)$, let $\rdes(\sigma)$  denote the vector difference $(\alpha_1,\dotsc,\alpha_k) - (1^k)$ (leaving in place any zeros created in the process). For example, $\rdes(132)=010$ and $\rdes(1432) = 0210$. Up to a relabelling, the weak compositions $\rdes(\sigma)$ are the ones introduced by Garsia in \cite{Gar:1980}. They are enumerated by $[n]_q!$.

Bergeron and Reutenauer define their sets $\PI_n$ recursively in such a way that
\begin{itemize}
\item $\PI_0 := \{0\}$,
\item $1^n{+}\PI_{n-1} \subseteq \Destd_{(n)}$ and $\Destd_{(n)}$ is disjoint from $\PI_{n-1}$, and
\item $\PI_n := \PI_{n-1} \cup \Destd_{(n)} \setminus \bigl(1^n{+}\PI_{n-1}\bigr) .$
\end{itemize}
Here, $1^n{+}\PI_{n-1}$ denotes the vector sums $\{(1^n) + \destd \suchthat \destd \in \PI_{n-1}\}$. Note that the compositions in $\Destd_{(n)}$ all have length $n$. Moreover, $1^{n+1}{+}\Destd_{(n)} \subseteq \Destd_{(n+1)}$. Indeed, if $\sigma=\sigma'1$ is a permutation in $\frakS_{n+1}$ with suffix ``1'' in one-line notation, then $(1^{n+1}) + \destd(\std(\sigma')) = \destd(\sigma)$. That \eqref{eq: recursion} enumerates $\PI_n$ is immediate \cite[Proposition 6.1]{GarWal:2003}. We give the first few sets $\PI_n$ and $\Destd_{(n)}$ in Figure \ref{fig: example-comps}.
%


\begin{figure}[hbt]
\vskip-.2in
{\small
\begin{align*}
   \Destd_{(1)} &= \{\underline{1}\} & \PI_0 &=\{0\} \\
   \Destd_{(2)} &= \{\underline{11},21\}  & \PI_1 &= \{0\} \\
   \Destd_{(3)} &= \{\underline{111}, 211, 121, 221, 212, \underline{321}\} & \PI_2 &= \{0,21\} \\
   \Destd_{(4)} &=
   \{\underline{1111}, 2111, 1211, 1121, 2211, 2121, 1221, 2112, 1212, 2221, 2212, 2122, &
   \PI_3 &= \{0,21, 211, \\
   & \phantom{\ =1} \underline{3211}, 3121, 1321, \underline{3221}, \underline{2321}, 3212, 2312, 2132, \underline{3321}, \underline{3231}, 3213, 4321\}&&
   \phantom{\ =\ \ } 121,221, 212\}
\end{align*}
}
\vskip-.2in
\caption{The sets $\Destd_{(n)}$ and $\PI_n$ for small values of $n$.
Compositions $1^n{+}\PI_{n-1}$ are underlined in $\Destd_{(n)}$.}
\label{fig: example-comps}
\end{figure}

\subsection{Pure and inverting compositions}

We now give an alternative description of the compositions in $\PI_n$ introduced by Bergeron and Reutenauer which will be easier to work with in what follows.
Call a composition $\alpha$ {\it inverting} if and only if for each $i>1$ (with $i$ less than or equal to the largest part of $\alpha$) there exists a pair of indices $s<t$ such that $\alpha_s=i$ and $\alpha_t=i-1$.  For example, $13112312$ is inverting while $21123113$ is not.  Any composition $\alpha$ admits a unique factorization
\begin{gather}\label{eq: pure}
	\alpha = \gamma k^{i_k} \dotsb 2^{i_2}1^{i_1} \ \quad (i_j \ge 1) ,
\end{gather}
such that $\gamma$ is a composition that does not contain any of the values from $1$ to $k$, and $k$ is maximal (but possibly zero).  We say $\alpha$ is {\it pure} if and only if this maximal $k$ is even.  (Note that if the last part of a composition is not $1$, then $k=0$ and the composition is pure.)  For example, $5435211$ is pure with $k=2$ while $3231$ is impure since $k=1$.

\begin{proposition}[\cite{BerReu}]
The set of inverting compositions of length $n$ is precisely $\Destd_{(n)}$.
The set of pure and inverting compositions of length at most $n$ is precisely $\PI_n$.
\end{proposition}

We reprise the proof of Bergeron and Reutenauer, for the sake of completeness.

\begin{proof}
Let $\mathcal D_{(n)}$ denote the set of inverting compositions of length $n$. The destandardization procedure makes it clear that $\Destd_{(n)} \subseteq \mathcal D_{(n)}.$ For the reverse containment, we use induction on $n$ to show that $|\mathcal D_{(n)}|=n!$. (The base case $n=1$ is trivially satisfied.) Let $\alpha=(a_1,\dotsc,a_{n-1})$ be one of the $(n-1)!$ compositions in $\mathcal D_{(n-1)}$. We construct $n$ distinct compositions by inserting a new part between positions $k$ and $k+1$ in $\alpha$ (for all $0\leq k <n-1$). Define this part $m_k(\alpha)$ by
\[
	m_k(\alpha) = \max\bigl(\{a_i \suchthat i\leq k\} \cup \{1+a_j \suchthat j>k\}\bigr).
\]
To reverse the procedure, simply remove the rightmost maximal value appearing in the inverting composition of length $n$. Conclude that applying the procedure to $\mathcal D_{(n-1)}$ results in $n!$ distinct elements in $\mathcal D_{(n)}$. Finally, since the reverse map from $\mathcal D_{(n)}$ to $\mathcal D_{(n-1)}$ is an $n$ to $1$ map, we get that $|\mathcal D_{(n)}| = n!$.

Turning to $\PI_n$, we argue  that $\PI_n \cap \Destd_{(n)}$ are the pure compositions in $\Destd_{(n)}$ of length $n\geq0$.
This will complete the proof, since by construction and the previous paragraph, the compositions $\PI_n$ are inverting. (Indeed, $\PI_{n}\subseteq \bigcup_{0\leq i\leq n} \Destd_{(i)}$, setting $\Destd_{(0)}=\{0\}$.) We argue by induction on $n$. (The base case $n=0$ is trivially satisfied.) Note that if $\alpha \in \Destd_{(n)}$ is impure, then $k$ is odd in the factorization \eqref{eq: pure}, and $\alpha':=\alpha - (1^n)$ is pure. That is, $\alpha' \in \PI_{n-1}\subseteq \PI_n$. These are precisely the compositions eliminated from $\Destd_{(n)}$ in constructing $\PI_n$, for $\PI_n := \PI_{n-1} \cup \Destd_{(n)} \setminus \bigl(1^n{+}\PI_{n-1}\bigr)$. In other words, if $\alpha \in \Destd_{(n)}$ is pure, then $\alpha\in\PI_n$.
\end{proof}

\subsection{A bijection}

Let $\Comp_{n,d}$ be the set of all compositions of $d$ into at most $n$ parts and set $\PPI_{n,d}:= \{ (\lambda, \beta) \suchthat \lambda \hbox{ a partition}, \beta \in \PI_n, |\lambda| + |\beta|=d, \hbox{ and } \length(\lambda) \le n, \length(\beta) \le n \}$.  We define a map $\phi\colon \PPI_{n,d} \rightarrow \Comp_{n,d}$ as follows.

Let $(\lambda, \beta)$ be an arbitrary element of $\PPI_{n,d}$.  Then $\phi((\lambda, \beta))$ is the composition obtained by adding $\lambda_i$ to the $i^{\textit{th}}$ largest part of $\beta$ for each $1 \le i \le \length(\lambda)$, where if $\beta_j=\beta_k$ and $j<k$, then $\beta_j$ is considered smaller than $\beta_k$.  If $\length(\lambda) > \length(\beta)$, append zeros after the last part to lengthen $\beta$ before applying $\phi$.  (See Figure \ref{leadingterm}.)

\begin{figure}[hbt]
\centering{\small
\begin{tabular}{r c@{\ \ }c@{\ \ }c@{\ \ }c@{\ \ }c@{\ \ }c@{\ \ }c@{\ \ }c@{\ \ }c@{\ \ }c@{\ \ }c}
$\lambda\ =\ $ & 1 & 4 & 2 & 1 & 1 & 4 & 5 & 2 & 4 & 1 & 1 \\
$\beta\ =\ $ &  2 & 4 & 3 & 1 & 1 & 3 & 4 & 2 & 3 &  &   \\
\hline
$\phi(\lambda, \beta) \ =\ $ & 3 & 8 & 5 & 2 & 2 & 7 & 9 & 4 & 7 & 1 & 1  \end{tabular}
}
\caption{An example of the map $\phi\colon \PPI_{13,49} \rightarrow \Comp_{13,49}$.}
\label{leadingterm}
\end{figure}

\begin{proposition}\label{thm: bijection}
The map $\phi$ is a bijection between $\PPI_{n,d}$ and $\Comp_{n,d}$
\end{proposition}

\begin{proof}
We prove this by describing the inverse $\phi^{-1}$ algorithmically.  Let $\alpha$ be an arbitrary composition in $\Comp_{n,d}$ and set $(\lambda, \beta) := (\emptyset, \alpha)$.
\begin{enumerate}
\item If $\beta$ is pure and inverting, then $\phi^{-1}(\alpha) := (\lambda,\beta)$
\item If $\beta$ is impure and inverting, then set $\phi^{-1}(\alpha) := (\lambda+(1^n), \beta-(1^n))$.
\item \label{itm: not-inverting} If $\beta$ is not inverting, then let $j$ be the smallest part of $\beta$ such that there does not exist a pair of indices $s<t$ such that $\beta_s=j$ and $\beta_t=j-1$.  Let $m$ be the number of parts of $\beta$ which are greater than or equal to $j$.  Replace $\beta$ with the composition obtained by subtracting $1$ from each part greater than or equal to $j$ and replace $\lambda$ with the partition obtained by adding $1$ to each of the first $m$ parts.
\item Repeat Steps (1)--(4) until $\phi^{-1}$ is obtained, that is, until Step (1) or (2) above is followed.
\end{enumerate}

To see that $\phi \phi^{-1} = \id$, consider an arbitrary composition $\alpha$.  If $\alpha$ is pure and inverting, then $\phi \phi^{-1}(\alpha)=\phi( \emptyset, \alpha) = \alpha$.  If $\alpha$ is impure and inverting, then $\phi(\phi^{-1}(\alpha))=\phi(((1^{\length(\alpha)}), \alpha-(1^{\length(\alpha)}) ))=\alpha$.  Finally, consider a composition $\alpha$ which is not inverting.  Note that the largest entry in $\alpha$ is decreased at each iteration of Step (3).  Therefore the largest entry in the partition records the number of times the largest entry in $\alpha$ is decreased.  Similarly, for each $i \le \length(\lambda)$, the $i^{\textit{th}}$ largest entry in $\alpha$ is decreased by one $\lambda_i$ times.  This means that the $i^{\textit{th}}$ largest part of $\alpha$ is obtained by adding $\lambda_i$ to the $i^{\textit{th}}$ largest part of $\beta$ and therefore our procedure $\phi^{-1}$ inverts the map $\phi$.
\end{proof}

Figure \ref{fig: phi-inverse} illustrates the algorithmic description of $\phi^{-1}$ as introduced in the proof of Proposition \ref{thm: bijection} on $\alpha=38522794711$.
\begin{figure}[hbt]
\vskip-.2in
{
\begin{align*}
\alpha\mapsto\begin{pmatrix}\lambda \\ \beta \end{pmatrix}:\ \  &
\begin{tabular}{c@{\;}c@{\;}c@{\;}c@{\;}c@{\;}c@{\;}c@{\;}c@{\;}c@{\;}c@{\;}c}
&&&&&$\emptyset$\\
3 & \underline{8} & \underline{5} & 2 & 2 & \underline{7} & \underline{9} & \underline{\underline{4}} & \underline{7} & 1 & 1
 \end{tabular}
\rightarrow
\begin{tabular}{c@{\;}c@{\;}c@{\;}c@{\;}c@{\;}c@{\;}c@{\;}c@{\;}c@{\;}c@{\;}c}
 & 1 & 1 & & & 1 & 1 & 1 & 1   \\
3 & \underline{7} & 4 &2 & 2 & \underline{\underline{6}} & \underline{8} & 3 & \underline{\underline{6}} & 1 & 1 \\[.5ex]
 \end{tabular}
\\
& \hskip2.25in \downarrow \\[.5ex]
&
\begin{tabular}{c@{\;}c@{\;}c@{\;}c@{\;}c@{\;}c@{\;}c@{\;}c@{\;}c@{\;}c@{\;}c}
  & 3 & 1 & & & 3 & 3 & 1 & 3   \\
3 & 5 & 4 & 2 & 2 & 4 & \underline{\underline{6}} & 3 & 4 & 1 & 1 \\[.5ex]
 \end{tabular}
\leftarrow
\begin{tabular}{c@{\;}c@{\;}c@{\;}c@{\;}c@{\;}c@{\;}c@{\;}c@{\;}c@{\;}c@{\;}c}
 & 2 & 1 & & & 2 & 2 & 1 & 2  \\
3 & \underline{6} & 4 & 2  & 2 & \underline{\underline{5}} & \underline{7} & 3 & \underline{\underline{5}} & 1 & 1 \\[.5ex]
 \end{tabular}
\\
& \hskip.65in {\downarrow} \\[.5ex]
&
\begin{tabular}{c@{\;}c@{\;}c@{\;}c@{\;}c@{\;}c@{\;}c@{\;}c@{\;}c@{\;}c@{\;}c}
   & 3 & 1 & & & 3 & 4 & 1 & 3 \\
3 & 5 & 4 & 2 & 2 & 4 & 5 & 3 & 4 & 1 & 1 \\[.5ex]
 \end{tabular}
\rightarrow
\begin{tabular}{c@{\;}c@{\;}c@{\;}c@{\;}c@{\;}c@{\;}c@{\;}c@{\;}c@{\;}c@{\;}c}
1 & 4 & 2 & 1 & 1 & 4 & 5 & 2 & 4 & 1 & 1  \\
2 & 4 & 3 & 1 & 1 & 3 & 4 & 2 & 3
 \end{tabular}
\rightarrow\ \begin{pmatrix} 54442211111 \\ 243113423 \end{pmatrix}.
\end{align*}
}
\vskip-.2in
\caption{The map $\phi^{-1}\colon \Comp_{13,49} \to \PPI_{13,49}$ applied to $\alpha=38522794711$. Parts $j$ from Step \ref{itm: not-inverting} of the algorithm are marked with a double underscore.}
\label{fig: phi-inverse}
\end{figure}

\section{Main Theorem}\label{sec: main-results}

Let $\PI_n$ be as in Section \ref{sec: coinvariant-space} and set $\cB_n := \{ \qschur_{\beta}\suchthat \beta \in \PI_n \}$. We prove the following.

\begin{theorem}\label{thm: main}
The set $\cB_n$ is a basis for the $Sym_n$-module $R_n$.
\end{theorem}

To prove this, we analyze the quasisymmetric polynomials $QSym_{n,d}$ in $n$ variables of homogeneous degree $d$. Note that ${QSym_n = \bigoplus_{d\geq0} QSym_{n,d}}$. Therefore, if $\sS_{n,d}$ is a basis for $QSym_{n,d}$, then the collection ${\bigcup_{d\geq0} \sS_{n,d}}$ is a basis for $QSym_n$. First, we introduce a useful term order.

\subsection{The revlex order}

Each composition $\alpha$ can be rearranged to form a partition $\shape(\alpha)$ by arranging the parts in weakly decreasing order.  Recall the {\it lexicographic} order $\lexgeq$ on partitions of $n$, which states that $\lambda \lexgeq \mu$ if and only if the first nonzero entry in $\lambda - \mu$ is positive.  For two compositions $\alpha$ and $\gamma$ of $n$, we say that $\alpha$ is larger then $\gamma$ in {\it revlex} order (written $\alpha \succeq \gamma$) if and only if either

\begin{itemize}
\item $\shape(\alpha) \lexgeq \shape(\gamma)$, or
\item $\shape(\alpha)=\shape(\gamma)$ and $\alpha$ is lexicographically larger than $\gamma$ when reading right to left.
\end{itemize}
For instance, we have
\[
   4 \succeq 13 \succeq 31 \succeq 22 \succeq 112 \succeq 121 \succeq 211 \succeq 1111.
\]

\noindent{\it Remark:\ }
Extend revlex to weak compositions of $n$ of length at most $n$ by padding the beginning of $\alpha$ or $\gamma$ with zeros as necessary, so $\length(\alpha)=\length(\gamma)=n$. Viewing these as exponent vectors for monomials in $\xx$ provides a term ordering on $\bQ[\xx]$. However, it is not good term ordering in the sense that it is not multiplicative: given exponent vectors ${\alpha}$, ${\beta}$, and ${\gamma}$ with ${\alpha} \succeq {\gamma}$, it is not necessarily the case that  ${\alpha}+{\beta} \succeq {\gamma}+{\beta}$. This is likely the  trouble encountered in \cite{BerReu} and \cite{GarWal:2003} when trying to prove the Bergeron--Reutenauer conjecture \eqref{conj:basis}. We circumvent this difficulty by working with the Schur polynomials $s_\lambda$ and the quasisymmetric Schur polynomials $\qschur_\alpha$. We consider leading polynomials $\qschur_\gamma$ instead of leading monomials $x^\gamma$. The leading term $\qschur_\gamma$ in a product $s_\lambda\cdot\qschur_\alpha$ is readily found.

\subsection{Proof of main theorem}

We claim that the collection $\sS_{n,d}=\{ s_{\lambda} \qschur_{\beta} \suchthat |\lambda|+|\beta|=d, \length(\lambda) \le n, \length(\beta) \le n, \hbox{ and } \beta \in \PI_n \}$ is a basis for $QSym_{n,d}$, which in turn implies that $\cB_n$ is a basis for $R_n$.  To prove this, we make use of a special \lrct called the {\it super filling}.  Consider a composition $\beta$ and a partition $\lambda=(\lambda_1, \hdots , \lambda_k)$.   If $\length(\lambda) > \length(\beta)$ then append $\length(\lambda)-\length(\beta)$ zeros to the end of $\beta$.  Fill the cells in the $i^{th}$ row from the bottom of $\beta$ with the entries $k+i$.  Append $\lambda_i$ cells to the $i^{th}$ longest row of $\beta$.  (If two rows of $\beta$ have equal length, the lower of the rows is considered longer.)  These new cells are then filled so that their entries have content $\lambda^*$ as follows.  Fill the new cells in the $j^{th}$ longest row with the entries $\lambda_{k-j+1}$ unless two rows are of the same length.  If two rows are the same length, fill the lower row with the lesser entries.  The resulting filling is called the {\it super filling $S(\lambda,\beta)$}.

\begin{proposition}{\label{prop: satisfies}}
The super filling $S(\lambda,\beta)$ obtained from composition $\beta$ and partition $\lambda$ is a filling satisfying \eqref{itm: LR-row}--\eqref{itm: LR-lattice}.
\end{proposition}

\begin{proof}
The super filling
$S(\lambda, \beta)$ satisfies \eqref{itm: LR-row} and \eqref{itm: LR-content} by construction.  We must prove that the filling also satisfies \eqref{itm: LR-CTs} and \eqref{itm: LR-lattice}.  Note that since $S( \lambda, \beta)$ satisfies \eqref{itm: CT-row} by construction, we need only prove that the entries in the filling satisfy the triple condition \eqref{itm: CT-triple} and the lattice condition \eqref{itm: LR-lattice}.  In the following, let $\alpha$ be the shape of $S(\lambda, \beta)$.

To prove that the filling $S(\lambda, \beta)$ satisfies \eqref{itm: CT-triple}, consider an arbitrary pair of cells $(i,k)$ and $(j,k)$ in the same column. If $\alpha_i \ge \alpha_j$ then $\beta_i \ge \beta_j$, since the entries from $\lambda$ are appended to the rows of $\beta$ from largest row to smallest row.  Therefore if $(i,k)$ is a cell in the diagram of $\beta$ then $T(j,k) < T(i,k)=T(i,k-1)$ regardless of whether or not $(j,k)$ is in the diagram of $\beta$.  If $(i,k)$ is not in the diagram of $\beta$ then $(j,k)$ cannot be in the diagram of $\beta$ since $\beta_i \ge \beta_j$.  Therefore $T(j,k) < T(i,k) $ since the smaller entry is placed into the shorter row, or the lower row if the rows have equal length.

If $\alpha_i < \alpha_j$ then $\beta_i  \le \beta_j$.  If $T(i,k) \le T(j,k)$ then $(i,k)$ is not in the diagram of $\beta$.  If $(j,k+1)$ is in the diagram of $\beta$ then $T(i,k) < T(j,k+1)$ since the entries in the diagram of $\beta$ are larger than the appended entries.  Otherwise the cell $(j,k+1)$ is filled with a larger entry than $(i,k)$ since the longer rows are filled with larger entries and $\alpha_j>\alpha_i$.  Therefore the entries in $S(\lambda,\beta)$ satisfy \eqref{itm: CT-triple}.

To see that the entries in $S(\lambda, \beta)$ satisfy \eqref{itm: LR-lattice}, consider an entry $i$.  We must show that an arbitrary prefix of the reading word contains at least as many $i$'s as $(i-1)$'s.  (Note that this is true when the prefix chosen is the entire reading word since $\lambda_i^* \ge \lambda_{i-1}^*$.)  Let $c_i$ be the rightmost column of $S(\lambda,\beta)$ containing the letter $i$ and let $c_{i-1}$ be the rightmost column of $S(\lambda,\beta)$ containing the letter $i-1$.  Note that all entries not in the diagram of $\beta$ in a given row are equal.  If $c_i >c_{i-1}$ then every prefix will contain at least as many $i's$ as $(i-1)$'s since there will always be at least one $i$ appearing before any pairs $i,i-1$ in reading order.  If $c_i=c_{i-1}$, then the entry $i$ will appear in a higher row than the entry $i-1$ and hence will be read first for each column containing both an $i$ and an $i-1$.  Therefore the reading word is a reverse lattice word and hence the filling satisfies \eqref{itm: LR-lattice}.
\end{proof}

\begin{proof}[Proof of Theorem \ref{thm: main}]

Order the compositions of $d$ into at most $n$ parts by the revlex order.  To define the ordering on the elements of $\sS_{n,d}$, note that their indices are pairs of the form $(\lambda, \beta)$, where $\lambda$ is a partition of some $k \le d$ and $\beta$ is a composition of $d-k$ which lies in $\PI_n$.  We claim that the leading term in the quasisymmetric Schur polynomial expansion of $s_{\lambda} \qschur_{\beta}$ is the polynomial $\qschur_{\phi(\lambda, \beta)}$.  To see this, recall from Proposition \ref{thm:LR-qschur} that the terms of $s_{\lambda} \qschur_{\beta}$ are given by \lrctx of shape $\alpha\supseteq\beta$ and appended content $\lambda^*$, where $\alpha$ is an arbitrary composition shape obtained by appending $|\lambda|$ cells to the diagram of $\beta$ so that conditions \eqref{itm: CT-row} and \eqref{itm: CT-triple} are satisfied.

To form the largest possible composition (in revlex order), one must first append as many cells as possible to the longest row of $\beta$, where again the lower of two equal rows is considered longer.  The filling of this new longest row must end in an $L:=\length(\lambda)$, since the reading word of the \lrct must satisfy \eqref{itm: LR-lattice}.  No entry smaller than $L$ can appear to the left of $L$ in this row, since the row entries are weakly decreasing from left to right.  This implies that the maximum possible number of entries that could be added to the longest row of $\beta$ is $\lambda_1$.  Similarly, the maximum possible number of entries that can be added to the second longest row of $\beta$ is $\lambda_2$ and so on.  If $\length(\lambda) > \length(\beta)$, append the extra parts of $\lambda$ (from least to greatest, top to bottom) after the bottom row of $\beta$.  The resulting shape is precisely the shape of $S(\lambda,\beta)$ which is equal to $\phi(\lambda, \beta)$ since $\beta$ is a pure and inverting composition.  Therefore there is at least one \lrct of the $\phi(\lambda, \beta)$ since $S(\lambda, \beta)$ is an \lrct by Proposition \ref{prop: satisfies}.

The shape of the \lrct $S(\lambda,\beta)$ corresponds to the largest composition appearing as an index of a quasisymmetric Schur polynomial in the expansion of $s_{\lambda} \qschur_{\beta}$, implying that $\qschur_{\phi(\lambda,\beta)}$ is indeed the leading term in this expansion.  Since $\phi$ is a bijection, the entries in $\sS_{n,d}$ span $\Qsym_{n,d}$ and are linearly independent.  Therefore $\sS_{n,d}$ is a basis for $\Qsym_{n,d}$ and hence $\cB_n$ is a basis for the $\Sym_n$-module $R_n$.
\end{proof}

\begin{remark}{\label{onlyfilling}}
{\rm%
Note that in the proof of Theorem \ref{thm: main}, the entries appearing in the filling of shape $\phi(\lambda, \alpha)$ are uniquely determined by the lattice condition \eqref{itm: LR-lattice}.  This implies that $C_{\lambda,\alpha}^{\phi(\lambda,\alpha)} = 1$.  This fact allows us to work over $\mathbb{Z}$, a slightly more general setting than working over $\mathbb{Q}$.  (See Section \ref{sec: integrality} for details.)
}
\end{remark}

The transition matrix between the basis $\sS_{3,4}$ and the quasisymmetric Schur polynomial basis for $QSym_{3,4}$ is given in Figure \ref{fig: example}.

\def\o{\hbox{\hskip.6em}\cdot\hbox{\hskip.6em}}
\begin{figure}[hbt]
 \centering
\begin{tabular}{r l}
& \begin{tabular}{@{\hskip1.55em} c @{\hskip2.15em}c @{\hskip2.05em} %
c @{\hskip1.85em}c @{\hskip1.6em}c @{\hskip1.5em}c @{\hskip1.4em} c }
 	4 & 13 & 31 & 22 & 112 & 121 & 211
  \end{tabular}\\[2ex]
\begin{tabular}{r}
 $ s_4$ \\[.2ex]
 $ s_{31}$ \\[.2ex]
 $s_1 {\,\cdot\,} \qschur_{21}$ \\[.2ex]
 $ s_{22}$ \\[.2ex]
 $ s_{211}$ \\[.2ex]
 $\qschur_{121} $ \\[.2ex]
 $\qschur_{211} $
\end{tabular} &
$\displaystyle\left(
\begin{array}{@{}ccccccc @{}}
 1 & \o & \o & \o & \o & \o & \o \\[.2ex]
 \o & 1 & 1 & \o & \o & \o & \o \\[.2ex]
 \o & \o & 1 & 1 & \o & \o & 1 \\[.2ex]
 \o & \o & \o & 1 & \o & \o & \o \\[.2ex]
 \o & \o & \o & \o & 1 & \o & \o \\[.2ex]
 \o & \o & \o & \o & \o & 1 & \o \\[.2ex]
 \o & \o & \o & \o & \o & \o & 1
\end{array}\right)$
\end{tabular}
\caption{The transition matrix for $n=3, d=4$.}
\label{fig: example}
\end{figure}

\section{Corollaries and applications}\label{sec: corollaries}

\subsection{Closing the Bergeron--Reuteuaner conjecture}
The relationship between the monomial basis and quasisymmetric Schur basis was investigated in \cite[Thm. 6.1 \& Prop. 6.7]{HLMvW:1}. We recall the pertinent facts.

\begin{proposition}[\cite{HLMvW:1}]\label{thm: qs-as-m}
The polynomials $M_\gamma$ are related to the polynomials $\qschur_\alpha$ as follows:
\begin{equation}\label{eq: S=M sum}
   \qschur_\alpha = \sum_\gamma K_{\alpha,\gamma}\,M_\gamma \,,
\end{equation}
where $K_{\alpha,\gamma}$ counts the number of composition tableaux $T$ of shape $\alpha$ and content $\gamma$. Moreover, $K_{\alpha,\alpha}=1$ and $K_{\alpha,\gamma}=0$ whenever $\shape(\alpha) \lexl \shape(\gamma)$.
\end{proposition}

We need a bit more to prove Conjecture \eqref{conj:basis}.

\begin{lemma}\label{thm: qs-as-m extension}
In the notation of Proposition \ref{thm: qs-as-m}, $K_{\alpha,\gamma} = 0$ whenever $\shape(\alpha)=\shape(\gamma)$ and $\alpha\neq \gamma$.
\end{lemma}

\begin{proof}
We argue by induction on the largest part of $\alpha$ that if $\shape(\alpha)=\shape(\gamma)$, and $T$ is a composition tableau with shape $\alpha$ and content $\gamma$, then $\alpha=\gamma$.

The base case is trivial, for if the largest part of $\alpha$ is $1$, then $\alpha=\gamma=(1^d)$ for some $d$. Now suppose $\alpha$ has largest part $l$. We claim that all rows $i$ in $T$ of length $l$ must be filled only with $i$'s. This claim finishes the proof. Indeed, we learn that $\alpha_i = \gamma_i$ for all such $i$. Thus we may apply the induction hypothesis to the new compositions $\alpha'$ and $\gamma'$ obtained by deleting the largest parts from each.


To prove the claim, suppose row $i$ of $T$ has length $l$ and is not filled with all $i$'s.  Let $(i,k)$ be the rightmost cell in row $i$ containing the entry $i$.  The $i$ in column $k+1$ must appear in a lower row, say row $j$, by condition \eqref{itm: CT-row} since the entries above row $i$ in the first column must be less than $i$.  This implies that $T(i,k)=T(j,k+1)$.  But $T(j,k) \ge T(j,k+1)$ and hence $T(j,k) \ge T(i,k)$, so \eqref{itm: CT-triple} is violated regardless of which row is longer.  Therefore row $i$ must be filled only with $i$'s and the claim follows by induction.
\end{proof}

\begin{theorem} \label{thm: m-as-qs}
In the expansion $M_{\alpha} = \sum_{\gamma} \tilde{K}_{\alpha, \gamma}\, \qschur_{\gamma},$
$\tilde{K}_{\alpha, \alpha} = 1$ and $\tilde{K}_{\alpha, \gamma} = 0$ whenever $\alpha \prec \gamma$.
\end{theorem}

\begin{proof}
From Proposition \ref{thm: qs-as-m} and Lemma \ref{thm: qs-as-m extension}, we learn that $K_{\alpha,\gamma}=0$ whenever $\alpha\prec\gamma$. (The proposition handles the first condition in the definition of the revlex order and the lemma handles the second condition.) Now arrange the integers $K_{\alpha,\gamma}$ in a matrix $K$, ordering the rows and columns by $\succeq$. The previous observation shows that this change of basis matrix is upper-unitriangular. Consequently, the same holds true for $\tilde K = K^{-1}$.
\end{proof}

We are ready to prove Conjecture \eqref{conj:basis}. Let $\PI_n$ and $R_n$ be as in Section \ref{sec: main-results}.

\begin{corollary}\label{thm: m-basis}
The set $\{ M_{\beta} \suchthat \beta \in \PI_n \}$ is a basis for the $Sym_n$-module $R_n$.
\end{corollary}

\begin{proof}
We show that the collection $\sM_{n,d}=\{ s_{\lambda} M_{\beta} \suchthat |\lambda|+|\beta|=d, \length(\lambda) \le n, \length(\beta) \le n, \hbox{ and } \beta \in \PI_n \}$ is a basis for $QSym_{n,d}$, which in turn implies that $\{ M_{\beta} \suchthat \beta \in \PI_n \}$ is a basis for $R_n$.  We first claim that the leading term in the quasisymmetric Schur polynomial expansion of $s_{\lambda} M_{\beta}$ is indexed by the composition $\phi(\lambda, \beta)$. The corollary will easily follow.

Applying Theorem \ref{thm: m-as-qs}, we may write $s_\lambda M_\beta$ as
\[
   s_{\lambda} M_{\beta}  =
   s_{\lambda} \qschur_{\beta} + \sum_{\beta \succ \gamma} \tilde{K}_{\beta,\gamma}\, s_{\lambda} \qschur_{\gamma} \,.
\]
Note that for any composition $\gamma$, the leading term of $s_\lambda S_\gamma$ is indexed by $\phi(\lambda,\gamma)$. This follows by the same reasoning used in the proof of Theorem \ref{thm: main}. To prove the claim, it suffices to show that $\beta \succ \gamma \implies \phi(\lambda, \beta) \succ \phi(\lambda, \gamma)$.

Assume first that $\shape(\beta)=\shape(\gamma)$.  Let $i$ be the greatest integer such that $\beta_i > \gamma_i$.  The map $\phi$ adds $\lambda_j$ cells to $\beta_i$ and $\lambda_k$ cells to $\gamma_i$, where $\lambda_j \ge \lambda_k$.  Therefore $\beta_i+\lambda_j > \gamma_i + \lambda_k$.  Since the parts of $\phi(\lambda, \beta)$ and $\phi(\lambda, \gamma)$ are equal after part $i$, we have $\phi(\lambda, \beta) \succeq \phi(\lambda, \gamma)$.

Next assume that $\shape(\beta) \succ \shape(\gamma)$.  Consider the smallest $i$ such that the $i^{\textit{th}}$ largest part $\beta_j$ of $\beta$ is not equal to the $i^{\textit{th}}$ largest part $\gamma_k$ of $\gamma$.  The map $\phi$ adds $\lambda_i$ cells to $\beta_j$ and to $\gamma_k$, so that $\beta_j + \lambda_i > \gamma_k + \lambda_i$.  Since the largest $i-1$ parts of $\phi(\lambda, \beta)$ and $\phi(\lambda, \gamma)$ are equal, we have $\shape(\phi(\lambda, \beta)) \succ \shape(\phi(\lambda, \gamma))$.

We now use the claim to complete the proof. Following the proof of Theorem \ref{thm: main}, we arrange the products $s_\lambda M_\beta$ as row vectors written in the basis of quasisymmetric Schur polynomials. The claim shows that the corresponding matrix is upper-unitriangular. Thus $\sM_{n,d}$ forms a basis for $QSym_{n,d}$,
as desired.
\end{proof}

\subsection{Triangularity}
It was shown in Section \ref{sec: main-results} that the transition matrix between the bases $\sS$ and $\left\{\qschur_\alpha\right\}$ is triangular with respect to the revlex ordering. Here, we show that a stronger condition holds: it is triangular with respect to a natural partial ordering on compositions.
Every composition $\alpha$ has a corresponding partition $\shape(\alpha)$ obtained by arranging the parts of $\alpha$ in weakly decreasing order.  A partition $\lambda$ is said to {\it dominate} a partition $\mu$ iff $\sum_{i=1}^k \lambda_i \ge \sum_{i=1}^k \mu_i$ for all $k$.  Let $C_{\lambda,\beta}^{\alpha}$ be the coefficient of $\mathcal{S}_{\alpha}$ in the expansion of the product $s_{\lambda}\,\mathcal{S}_{\beta}$.

\begin{theorem}
If $\shape(\alpha)$ is not dominated by $\shape(\phi(\lambda,\beta))$, then $C_{\lambda,\beta}^{\alpha}=0$.
\end{theorem}

\begin{proof}
Let $(\lambda, \beta)$ be an arbitrary element of $\PPI_{n,d}$ and let $\alpha$ be an arbitrary element of $\Comp_{n,d}$.  Set $\gamma:=\phi(\lambda, \beta)$.  If $\gamma \preceq \alpha$ then $C_{\lambda,\beta}^{\alpha}=0$ (by the proof of Theorem \ref{thm: main}) and we are done. 

Hence, assume that $\alpha \succ \phi(\lambda, \beta)=\gamma$ and that $\shape(\alpha)$ is not dominated by $\shape(\gamma)$.
Let $k$ be the smallest positive integer such that $\sum_{i=1}^k \shape(\alpha)_i > \sum_{i=1}^k \shape(\gamma)_i$.  (Such an integer exists since $\shape(\alpha)$ is not dominated by $\shape(\gamma)$.)  Therefore $\sum_{i=1}^k \shape(\alpha)_i - \sum_{i=1}^k \shape(\beta)_i > \sum_{i=1}^k \shape(\gamma)_i - \sum_{i=1}^k \shape(\beta)_i$ and there are more entries in the longest $k$ rows of $\alpha\supseteq\beta$ then there are in the longest $k$ rows of $\gamma\supseteq\beta$.  This implies that there are more than $\sum_{i=1}^k \lambda_i$ entries from $\alpha\supseteq\beta$ contained in the longest $k$ rows of $\alpha$, since there are $\sum_{i=1}^k \lambda_i$ entries in the longest $k$ rows of $\gamma \supseteq \beta$.  This implies that in a \lrct of shape $\alpha\supseteq\beta$, the longest $k$ rows must contain an entry less than $L-k+1$ where $L=\length(\lambda)$.

The rightmost entry in the $i^{\textit{th}}$ longest row of $\alpha\supseteq\beta$ must be $L-i+1$ for otherwise the filling would not satisfy the reverse lattice condition.  This means that the longest $k$ rows of $\alpha$ must contain only entries greater than or equal to $L-i+1$, which contradicts the assertion that an entry less than $L-k+1$ appears among the $k$ longest rows of $\alpha$.  Therefore there is no such \lrct of shape $\alpha$ and so $C_{\lambda,\beta}^{\alpha}=0$.
%
\end{proof}

\subsection{Integrality} {\label{sec: integrality}}
Up to this point, we have been working with the symmetric and quasisymmetric polynomials over the rational numbers, but their defining properties are equally valid over the integers. Briefly, bases for $\Sym_n(\bZ)$ and $\Qsym_n(\bZ)$ are the Schur polynomials $s_\lambda$ and the monomial quasisymmetric polynomials $M_\alpha$, respectively. See \cite{Mac:1995} and \cite{Haz:2003} for details. 

\begin{lemma} The polynomials $\{\cS_\alpha \suchthat \length(\alpha)\leq n \}$ form a basis of $\Qsym_n(\bZ)$. 
\end{lemma}

\begin{proof} 
Proposition \ref{thm: qs-as-m extension} states that the change of basis matrix $K$ from $\{\cS_\alpha\}$ to $\{M_\alpha\}$ is upper-unitriangular and integral. In particular, $K$ is invertible over $\bZ$, meaning that $\{\cS_\alpha\}$ is a basis for $\Qsym_n(\bZ)$.
\end{proof}

One consequence of the proof of Theorem \ref{thm: main} is that $C_{\lambda,\beta}^{\phi(\lambda,\beta)}=1$.  (See Remark \ref{onlyfilling}.)  We exploit this fact below to prove stronger versions of Conjectures \eqref{conj:free} and \eqref{conj:basis}. 

\begin{corollary}\label{thm: over Z}
The algebra $\Qsym_n(\bZ)$ is a free module over $\Sym_n(\bZ)$. A basis is given by $\{s_\lambda \cS_\beta \suchthat \beta\in\Pi_n, \length(\lambda)\leq n, \hbox{ and } \length(\beta)\leq n\}$. 
Replacing $\cS_\beta$ by $M_\beta$ results in an alternative basis. 
\end{corollary}

\begin{proof} 
Theorem \ref{thm: main} combines with Proposition \ref{thm:LR-qschur} (and the fact that $C_{\lambda,\beta}^{\phi(\lambda,\beta)}=1$) to establish an upper-unitriangular, integral change of basis matrix $C$ between $\{\cS_\alpha \suchthat \length(\alpha)\leq n\}$ and $\{s_\lambda \cS_\beta \suchthat \beta\in\Pi_n, \length(\lambda)\leq n, \hbox{ and } \length(\beta)\leq n\}$. Since the former is an integral basis for $\Qsym_n(\bZ)$, so is the latter. Composition of $K$, $C$ and $K^{-1}$ establishes the result for the monomial quasisymmetric polynomials.
\end{proof}

{\small
\bibliographystyle{abbrv}
\bibliography{bibl}
\label{sec:biblio}
}

\end{document}